\documentclass[aps,pra,twocolumn,showpacs,amsmath,amssymb,superscriptaddress]{revtex4}
\usepackage{amsthm,amsmath,amsfonts,amssymb,amsxtra,appendix,bookmark,dsfont,latexsym,stix}

\makeatletter
\usepackage{hyperref}
\usepackage{delarray,color}
\usepackage{euscript}

\theoremstyle{plain}

\theoremstyle{remark}

\def\geqslant{\ge}
\def\leqslant{\le}
\def\bq{\begin{eqnarray}}
	\def\eq{\end{eqnarray}}
\def\bqq{\begin{eqnarray*}}
	\def\eqq{\end{eqnarray*}}

\def\eps{\varepsilon}

\newcommand{\im}{\mathrm{i}}

\renewcommand{\epsilon}{\varepsilon}
\newcommand{\bx}{\mathbf{x}}

\newcommand{\by}{\mathbf{y}}

\def\R {\mathbb{R}}

\def\cE {\mathcal{E}}

\def\R {\mathbb{R}}
\def\C {\mathbb{C}}

\renewcommand{\leq}{\leqslant}
\renewcommand{\geq}{\geqslant}

\newcommand{\bA}{\mathbf{A}}
\newcommand{\bT}{\mathbf{T}}

\newcommand{\bJ}{\mathbf{J}}
\newcommand{\bj}{\mathbf{j}}

\newcommand{\curl}{\mathrm{curl}}
\newcommand{\sgn}{\mathrm{sgn}}
\newcommand{\gt}{\widetilde{g}}

\usepackage{color}

\date{January, 2026}

\begin{document} 

\title{The 2D Chern-Simons-Schr\"odinger system reduced to 1D}

\author{Nicolas Rougerie}
\affiliation{Ecole Normale Sup\'erieure de Lyon \& CNRS, UMPA (UMR 5669)}
\email{nicolas.rougerie@ens-lyon.fr}

\author{Qiyun Yang}
\affiliation{Ecole Normale Sup\'erieure de Lyon, UMPA (UMR 5669)}
\email{qiyun.yang@ens-lyon.fr}

\begin{abstract}
We study a mean-field model for a system of 2D abelian anyons, given by the dynamics of a Schr\"odinger matter field coupled to a Chern-Simons gauge field. We derive an effective 1D equation by adding a strongly anisotropic trapping potential (waveguide) acting on the Schr\"odinger field, and tracing out the tight confinement direction. The effective dynamics in the loose direction of the waveguide turns out to be governed by the classical 1D quintic NLS equation. 
\end{abstract}


\maketitle

Effective magnetic flux attachment~\cite{DalGerJuzOhb-11,EdmEtalOhb-13,ValWesOhb-20,LunRou-16,ZhaSreGemJai-14,ZhaSreJai-15} via the interaction with artificial gauge fields is one of the chief proposals to mimic, with cold atoms, physics such as that of the fractional quantum Hall effect, or anyon exchange statistics. A very natural classical field theory for such a mechanism in 2D is given in terms of a Schr\"odinger field $\Psi(t,\bx)\in \C$ 
coupled to a gauge vector potentiel $\bA (t,\bx)\in  \R^{2}$ in the manner 
\begin{equation}\label{eq:flux attach}
\curl_\bx \, \bA  = 2\pi \beta |\Psi| ^2
\end{equation}
with $\beta$ an effective coupling constant. In a mean-field approximation of 2D anyons~\cite{Wilczek-90,Myrheim-99} we typically have $\beta \propto \alpha N$ with $\alpha \in [0,2]$ the exchange statistics parameter and $N$ the particle number (and the convention that the field $\Psi$ is $L^2$-normalized). The most usual choice of gauge for $\bA$ is 
\begin{align}\label{eq:gauge2D}
\bA (\bx) = \beta \int_{\R^2} \frac{(\bx-\by) ^\perp}{|\bx-\by|^2} |\Psi(\by)|^2 d\by = \beta \nabla^{\perp} w_0 \star |\Psi|^2
\end{align}
where $w_0 (\bx) = - \log |\bx|$, $\star$ denotes convolution and $\nabla^\perp$ the gradient rotated by $\pi/2$.

The energy/Hamiltonian governing the dynamics is given by minimal coupling in the manner 
\begin{equation}\label{eq:Hamiltonian}
\cE^{\mathrm{2D}} [\Psi] = \int_{\R^2} \left|\left(-\im \nabla_\bx + \bA \right) \Psi \right| ^2 + \frac{g}{2} |\Psi|^4  
\end{equation}
where $g$ measures short-range two-particles interactions. The Schr\"odinger equation of motion is 
\begin{align}\label{eq:2DPDEpre}
    \mathrm{i}  \partial_t \Psi  & = \partial_{\overline{\Psi}} \mathcal{E}^{2\mathrm{D}} [\Psi] \nonumber \\
    & = \left[ \left( -\mathrm{i}\nabla +\bA \right)^2  + g|\Psi|^2 \Psi - 2\beta (\nabla^{\perp} w_0) \star \mathbf{J}  \right] \Psi
\end{align}
with the current 
\begin{equation}\label{eq:current}
    \mathbf{J}  = \frac{1}{2} \left[ \overline{\Psi} \left( -\mathrm{i} \nabla  + \bA \right) \Psi  +  {\Psi} \overline{\left( -\mathrm{i} \nabla  + \bA \right) \Psi } \right].
\end{equation}
Note that the above model is mass-critical: for attractive interactions $g<0$ it can become unstable. We refer to~\cite{Ataei-24,AtaGirLun-25} for a thorough discussions of regions of stability/instability. We henceforth  set $\int_{\R^2} |\Psi|^2 = 1$ and assume $g\geq 0$ or $|g|$ small enough, which implies stability.

The above arises as a mean-field approximation~\cite{ChiSen-92,CorDubLunRou-19,LunRou-15,Girardot-19,GirLee-24,AtaGirLun-25,Visconti-25} to the many-body problem for 2D abelian anyons, where one perturbs around the bosonic end (almost bosonic anyons) and is proposed as an effective theory of the fractional quantum Hall effect~\cite{LopFra-91,Zhang-92,ZhaHanKiv-89}. It is also a reformulation of the Jackiw-Pi model with (non-dimensionalized) Lagrangian density~\cite{JacPi-90b,Dunne-95} 
\begin{equation}\label{eq:Lagrangian}
\im \overline{\Psi} D_t \Psi  - \frac{1}{2} |D_\bx \Psi| ^2 - \frac{g}{2} |\Psi| ^4 + \frac{1}{2} \sum_{0 \leq a,b,c \leq 2} \epsilon^{abc} \bA_a \partial_b \bA_c
\end{equation}
Here $\bA$ also has a time component denoted $\bA_0$, $\epsilon^{abc}$ is the antisymmetric Levi-Civita symbol and the covariant derivatives are given as 
\begin{align}\label{eq:cov der}
D_t \Psi &= \partial_t + \im \bA_0\nonumber \\
D_\bx \Psi&= (D_x \Psi, D_y \Psi) = \left(\partial_x + \im \bA_1, \partial_y + \im \bA_2\right). 
\end{align}
Variations of the Lagragian lead to coupled equations for $\Psi$ and $\bA$ that one can, in Coulomb gauge, explictly solve for $\bA$, leading to~\eqref{eq:2DPDEpre}, see e.g.~\cite{Ataei-24,Ataei-25,AtaLunNgu-24,BerBouSau-95,BerBouSau-95a,LiuSmiTat-14,RajSig-20} and references therein.

In this note we investigate a dimensional reduction of the above dynamics, obtaining new approximate solutions. There are several distinct motivations for this: 

\medskip

\noindent (i) While the above model (or more precisely its second quantization) is the most accepted for 2D anyons, there are quite a few candidates for 1D anyons: Calogero or Lieb-Liniger models~\cite{LeiMyr-77,LeiMyr-88,Polychronakos-89,PosThoTra-17,Ouvry-07}, Kundu anyons/chiral BF model/anyon Hubbard model ~\cite{Rabello-96,AglGriJacPiSem-96,Kundu-99,BatGuaOel-06,GreSan-15,TanEggPel-15,Bon-etal-Pos-21} ... it makes sense~\cite{HanLeiMyr-92,Sen-94} to wonder if one is the limit of the 2D model. The answer turns out negative, see~\cite{RouYan-23a,RouYan-23b} and our findings below.

\medskip

\noindent(ii) Anyons of the FQHE are introduced as 2D bulk quasi-particles~\cite{AroSchWil-84,LunRou-16} but mostly probed via their propagation along 1D edge channels~\cite{BarEtalFev-20,NakEtalMan-20}. Even though the edge physics of a FQHE sample is a delicate issue, it seems natural to look for 1D restrictions of the 2D anyon bulk model.

\medskip

\noindent (iii) In putative cold atoms experiments, an extra trapping potential can be added to study 1D physics. This in particular lead to the realization of the Lieb-Liniger and Tonks-Girardeau models of the 1D Bose gas (see e.g.~\cite[Section~V]{BloDalZwe-08} and references therein).

\medskip

We follow the latter train of thought by adding an external static scalar potential acting on the Schr\"odinger field $\Psi$, tightly confining it along the $y$ direction. We derive the effective 1D theory obtained by tracing out the confined coordinate, which turns out to be a quintic NLS equation.

Regarding practical realizations of Chern-Simons-Schr\"odinger (CSS) theory with cold atoms, we mention two main proposals. (i) Realizing a Laughlin state with a majority species of atoms and couple it to a minority species, which, for strong inter-species interactions, will bind to Laughlin quasi-holes~\cite{LunRou-16,CooSim-15,ZhaSreGemJai-14}. Few-particles Laughlin states of cold atoms have been realized recently~\cite{LeoEtalGrei-23,LunEtalJoc-24}. (ii) Coupling a Bose gas of two-level atoms to laser fields resonant with the transition~\cite{ValWesOhb-20}. A particular design of the laser can produce, in the adiabatic regime, a density-dependent gauge potential akin to that of CSS theory. This is very close in spirit to proposals leading to the 1D chiral BF model~\cite{EdmValOhb-15,ChiEtalCel-22} that have been experimentally realized~\cite{FroEtalTar-22}. 

In both cases, a mean-field approximation (assuming Bose condensation) leads to CSS as the effective classical field theory. We next proceed to investigate its' reduction to 1D.

\bigskip 

\noindent\textbf{Model in a tight wave guide.} We add the potential 
$$V_\epsilon (\bx) = V_\epsilon (x,y)= x^2 + \frac{y^2}{\epsilon^2} - \frac{1}{\epsilon}$$
to~\eqref{eq:2DPDEpre}, leading to the equation of motion (with a corresponding modification of the energy~\eqref{eq:Hamiltonian})
\begin{equation}\label{eq:2DPDE}
    \mathrm{i}  \partial_t \Psi  = \left[ \left( -\mathrm{i}\nabla +\bA \right)^2 + V_\epsilon \Psi + g|\Psi|^2 \Psi - 2\beta (\nabla^{\perp} w_0) \star \mathbf{J}  \right] \Psi
\end{equation}
Here $\epsilon$ is thought of as a small parameter, and hence we consider tight trapping in the $y$ direction. We took an harmonic potential for convenience, but this has no influence on the following discussion. The one-particle Hamiltonian in the $y$ direction
\begin{equation}\label{eq:HO}
-\partial_y ^2 + \frac{y^2}{\eps ^2} - \frac{1}{\epsilon}
\end{equation}
has ground state energy $0$ and a gap above that, of the order 
$$e_\eps = \epsilon ^{-1}.$$ 
For low energy states, the motion will be frozen in the ground state 
$$u_\eps (y)= \left(\pi \eps\right) ^{-1/4} e^{-\frac{y^2}{2\eps}}.$$ 
To study this it is useful to change gauge, realizing~\eqref{eq:flux attach} in a manner more adapted to the 1D geometry~\cite{Sen-94,RouYan-23a,RouYan-23b}. We set 
\begin{equation}\nonumber
    \bT_0(x,y) = 
    \left( \begin{matrix}
        -\pi \operatorname{sgn}(y)\delta_{x=0} \\
        0
    \end{matrix}\right) 
\end{equation}
so that 
\begin{equation}\label{eq:gauge1D}
\bT :=  \beta \, \bT_0 \star |\Psi|^2
\end{equation}
satisfies
\begin{equation}
 \curl_\bx \, \bT= 2\pi \beta |\Psi|^2
\end{equation}
just as $\bA$ did. Explicitly, the change of gauge 
\begin{multline}\label{eq:changegauge}
\Psi \to \Psi e^{-\im \beta S} \\ \mbox{ with } S = S_0 \star |\Psi|^2, \quad  S_0(\bx) := \arctan \frac{y}{x} 
\end{multline}
allows us to work with $\bT$ as vector potential instead of $\bA$. Our ansatz for a quasi-1D solution will then be 
\begin{align}\label{eq:ansatz}
\Psi(t,x,y) &\simeq \varphi(t,x) u_{\varepsilon}(y) e^{-\mathrm{i}\beta S (t,x,y)}\nonumber\\
\quad S  &= S_0 \star |\varphi u_{\varepsilon}|^2 
\end{align}
and we shall solve for $\varphi (t,x)$ by projecting the 2D equation on the span of $u_{\varepsilon}$. We claim that the appropriate effective equation is the following cubic-quintic non-linear Schr\"odinger equation
\begin{equation}\label{eq:1DNLS}
 \boxed{\im \partial_t \varphi = - \partial_x ^2 \varphi + x^2 \varphi + \pi ^2 \beta ^2 |\varphi|^4 \varphi + \gt |\varphi|^2 \varphi}
\end{equation}
with $\beta$ the proportionality constant between magnetic field and matter density from~\eqref{eq:flux attach} and 
$$ \gt = g \int_{\R} |u_\eps| ^4.$$
A mathematically rigorous derivation in the $\eps \to 0$ limit is presented in~\cite{Yang-25}. In this note we present the key steps of the calculation. Then we proceed to investigate some static and dynamic properties of the 1D NLS~\eqref{eq:1DNLS}.

\bigskip 

\noindent\textbf{Effective energy.} Let us first compute the conserved energy~\eqref{eq:Hamiltonian} for our ansatz~\eqref{eq:ansatz}. Since 
\begin{equation}\nonumber
     \left( - \mathrm{i} \nabla_{\bx} +   \mathbf{A}\right) \Psi  =   e^{-\mathrm{i}\beta S} \left( - \mathrm{i} \nabla_{\bx} +  \mathbf{T} \right) (\varphi u_{\varepsilon })  ,
\end{equation}
and 
\begin{align*}
 \left( - \mathrm{i} \nabla_{\bx} +  \mathbf{T} \right)_x &= -\im \partial_x - \pi \beta \int_{\R} \sgn(y-y') |\Psi (x,y')|^2 \mathrm{d}y' \\
 \left( - \mathrm{i} \nabla_{\bx} +  \mathbf{T} \right)_y&= -\im \partial_y 
\end{align*}
a direct calculation gives 
\begin{align}\label{eq:calcener}
    \mathcal{E}^{2\mathrm{D}} [\Psi] &= \int_{\mathbb{R}^2 } \left| \left( - \mathrm{i} \nabla_{\bx} + \bT \right) (\varphi u_{\varepsilon })  \right|^2  \nonumber\\
    &+ \int_{\mathbb{R}^2} V_{\varepsilon} |\varphi u_{\varepsilon }|^2 + \frac{g}{2} \int_{\mathbb{R}^2} |\varphi u_\eps |^4  \nonumber\\
    &= \int_{\mathbb{R}} |x|^2 |\varphi(x)|^2\mathrm{d}x +\frac{\gt}{2} \int_{\mathbb{R}}|\varphi|^4  \nonumber\\ 
    &+ \int_{\mathbb{R}^2} \left| \left( \mathrm{i}\partial_x + \pi \beta f(y) |\varphi(x)|^2 \right) \varphi(x)\right|^2  u^2_{\varepsilon}(y) \mathrm{d}x \mathrm{d}y  
\end{align}
where we used the ground-state equation
\begin{equation}\label{eq:HOGS}
 - \partial_y ^2 u_\eps + \frac{y^2}{\eps ^2}u_\eps = \eps^{-1} u_\eps 
\end{equation}
and have denoted 
\begin{equation}\label{eq:f}
    f(y) =  \int_{\mathbb{R}} \sgn (y - y' ) u^2_{\varepsilon}(y') \mathrm{d} y'.
\end{equation}
Observe that 
\begin{equation}\label{eq:f2}
\partial_y f = 2 u_{\varepsilon}^2, \quad f( +\infty) = 1, \quad  f(-\infty) = -1
\end{equation}
so that direct integration gives
\begin{align}\label{eq:f3}
    \int_{\mathbb{R}} f(y) u^2 _{\varepsilon}(y)\mathrm{d} y &= 0 \nonumber \\
    \int_{\mathbb{R}} f^2(y)  u^2_{\varepsilon}(y) \mathrm{d} y &= \frac{1}{3}.
\end{align}
Expanding the square of the last term in~\eqref{eq:calcener} and inserting~\eqref{eq:f3} we find
\begin{equation}\nonumber
      \mathcal{E}^{2\mathrm{D}}[\Psi] = \mathcal{E}^{1\mathrm{D}} [\varphi] ,
\end{equation}
with
\begin{multline}\label{eq:1Dener}
    \mathcal{E}^{1\mathrm{D}} [\varphi] = \\\int_{\mathbb{R}} \left( |\partial_x \varphi (x)|^2  +  |x|^2 |\varphi(x)|^2 + \frac{1}{3} \pi^2 \beta^2 |\varphi(x)|^6  +\frac{\gt}{2} |\varphi|^4 \right) \mathrm{d} x.
\end{multline}
Hence, as long as 
\begin{equation}\label{eq:energybound}
\left|\mathcal{E}^{1\mathrm{D}} [\varphi] \right| \ll \eps^{-1}, 
\end{equation}
the gap of the transverse harmonic oscillator~\eqref{eq:HO}, the ansatz~\eqref{eq:ansatz} stays sensible, although not exact. And since~\eqref{eq:1Dener} is the conserved energy for~\eqref{eq:1DNLS} it is sufficient to have this bound for the initial datum. In~\cite{Yang-25} the second author gives a rigorous mathematical proof of the validity of the ansatz~\eqref{eq:ansatz} to approximate ground state configurations of our model. It is then natural to use this form as an initial datum for the dynamics, if the latter is studied at low enough energies.

\bigskip

\noindent\textbf{Effective dynamics.}  Henceforth we assume~\eqref{eq:energybound} and we proceed to the derivation of~\eqref{eq:1DNLS}. We multiply both sides of~\eqref{eq:2DPDE} by 
$${u_{\varepsilon}}e^{\mathrm{i} \beta S}$$
and then integrate over $y$. We discuss the different terms separately.

\medskip

\noindent$\bullet$ For the time derivative part, we have
\begin{equation}\nonumber
\begin{split}
    \int_{\mathbb{R}} & (\mathrm{i} \partial_t \Psi) {u_{\varepsilon}}e^{\mathrm{i} \beta S} \mathrm{d} y \\
    & =  \int_{\mathbb{R}} \mathrm{i} \Big[(\partial_t \varphi){u_{\varepsilon}}e^{ - \mathrm{i} \beta S}  \\
    &+ \left( - \mathrm{i} \beta \partial_t S\right) \varphi {u_{\varepsilon}}e^{ - \mathrm{i} \beta S} \Big]{u_{\varepsilon}}e^{\mathrm{i} \beta S} \mathrm{d} y \\
    & = \mathrm{i} \partial_t \varphi  + \beta  \varphi \int_{\mathbb{R} } S_0 \star (\partial_t (|\varphi u_{\varepsilon} |^2) ) u_{\varepsilon}^2 \mathrm{d} y\\
    & = \mathrm{i} \partial_t \varphi  + \beta  \varphi \int_{\mathbb{R} } S_0 \star (-2 \nabla \cdot \mathbf{J} ) u_{\varepsilon}^2 \mathrm{d} y \\
    & = \mathrm{i} \partial_t \varphi  - 2\beta  \varphi \int_{\mathbb{R} }\left( (\nabla S_0) \star \mathbf{J} \right) u_{\varepsilon}^2 \mathrm{d} y
\end{split}
\end{equation}
where we used the continuity equation 
$$ \partial_t |\Psi| ^2 + 2 \nabla_\bx \cdot \bJ = 0$$
derived from~\eqref{eq:2DPDE}, with $\bJ$ the current~\eqref{eq:current}.

\medskip

\noindent$\bullet$ For the current term (last term of the right-hand side of~\eqref{eq:2DPDE}) we have
\begin{multline}\nonumber
        \int_{\mathbb{R}} \left[  - 2\beta (\nabla^{\perp} w_0) \star \bJ  \right] {u_{\varepsilon}}e^{\mathrm{i} \beta S} \mathrm{d} y  \\ = - 2\beta  \varphi \int_{\mathbb{R}}   \left( (\nabla^{\perp} w_0) \star \bJ \right)   u_{\varepsilon}^2 \mathrm{d} y  \\
         = - 2\beta  \varphi \int_{\mathbb{R}}   \left( (\nabla S_0 +\mathbf{T}_0) \star \bj \right)   u_{\varepsilon}^2 \mathrm{d} y 
\end{multline}
where, with $\bT$ as in~\eqref{eq:gauge1D}, 
$$
\bj = \frac{1}{2} \left[ \overline{u_\eps \varphi} \left( -\mathrm{i} \nabla  + \bT \right) (u_\eps \varphi)  +  u_\eps \varphi \overline{\left( -\mathrm{i} \nabla  + \bT \right) (u_\eps \varphi) } \right].  
$$

\medskip

\noindent$\bullet$ For the potential part we directly have
\begin{multline*}
   \int_{\mathbb{R}}( V_{\varepsilon} \Psi) {u_{\varepsilon}}e^{ \mathrm{i} \beta S} \mathrm{d} y \\ = \int_{\mathbb{R}}\left( |x|^2 + \frac{1}{\varepsilon^2} |y|^2 \right) \varphi u_{\varepsilon}  e^{- \mathrm{i} \beta S} {u_{\varepsilon}} e^{\mathrm{i} \beta S} \mathrm{d} y\\
    = |x|^2 \varphi + \varphi \int_{\mathbb{R}} \frac{1}{\varepsilon^2}|y|^2 u_{\varepsilon}^2 \mathrm{d} y.
\end{multline*}

\medskip

\noindent$\bullet$ Similarly, for the $|\Psi|^4$ term
$$ g \int_{\R} |\Psi|^2 \Psi \, {u_{\varepsilon}}e^{\mathrm{i} \beta S} \mathrm{d} y = \gt |\varphi|^2 \varphi$$

\medskip

\noindent$\bullet$ As regards the magnetic kinetic operator we have 
\begin{align*}
        \int_{\mathbb{R} } & \left[ \left( -\mathrm{i}\nabla +\bA \right)^2 \Psi \right] {u_{\varepsilon}}e^{\mathrm{i} \beta S} \mathrm{d} y \\
        & = \int_{\mathbb{R} } \left[ \left( - \mathrm{i} \nabla_{} +  \bT\right)^2 (\varphi u_{\varepsilon })\right] {u_{\varepsilon}} \mathrm{d} y  \\
        & = \int_{\mathbb{R} } \Big[ \left( \mathrm{i}\partial_x + \pi \beta f(y) |\varphi(x,t)|^2 \right)^2 (\varphi(x,t) u_{\varepsilon }(y))\\  
        &- \partial^2_y  (\varphi(x,t) u_{\varepsilon }(y)) \Big] {u_{\varepsilon}(y)} \mathrm{d} y \\
        & = \int_{\mathbb{R} } \Big[ (-\partial_x^2 \varphi) u_{\varepsilon 
        } + \mathrm{i} \pi \beta f |\varphi|^2 (\partial_x \varphi) u_{\varepsilon} \\
        &+ \mathrm{i} \pi \beta f \partial_x (|\varphi|^2 \varphi) u_{\varepsilon} +  \pi^2 \beta^2 f^2 |\varphi|^4 \varphi u_{\varepsilon} + \varphi ( - \partial^2_y   u_{\varepsilon }) \Big] {u_{\varepsilon}} \mathrm{d} y  \\
        & = - \partial_x^2 \varphi + \frac{1}{3}\pi^2 \beta^2 |\varphi|^4 \varphi  + \varphi \int_{\mathbb{R} } ( - \partial^2_y   u_{\varepsilon }) u_{\varepsilon} \mathrm{d} y 
\end{align*}
where we used~\eqref{eq:f3}.

\medskip

Gathering the above computations and using~\eqref{eq:HOGS} we find
\begin{multline}\label{eq:preconc}
    \mathrm{i} \partial_t \varphi = - \partial_x^2 \varphi + x^2 \varphi + \frac{1}{3}\pi^2 \beta^2 |\varphi|^4  \varphi  + \gt |\varphi|^2 \varphi\\ - 2\beta  \varphi \int_{\mathbb{R}}   \bT_0 \star \bj  \, u_{\varepsilon}^2 \mathrm{d} y 
\end{multline}
and there remains to show that
\begin{equation}\label{eq:last claim} - 2\beta  \varphi \int_{\mathbb{R}}   \bT_0 \star \bj  \, u_{\varepsilon}^2 \mathrm{d} y =  \frac{2}{3}\pi^2 \beta^2 |\varphi|^4 \varphi. 
\end{equation}
Indeed, with $f$ as defined in~\eqref{eq:f}, we have
\begin{equation}\nonumber
   \bT =  \beta \left( \begin{matrix}
       -\pi f(y) |\varphi(x,t)|^2 \\
       0
   \end{matrix} \right)
\end{equation}
and hence
\begin{equation}\nonumber
    \begin{split}
        \bj & = \frac{1}{2} \left[ \overline{(\varphi u_{\varepsilon})} \left( -\mathrm{i} \nabla  + \bT) \right) (\varphi u_{\varepsilon})  +  {(\varphi u_{\varepsilon})} \overline{\left( -\mathrm{i} \nabla  + \bT) \right) (\varphi u_{\varepsilon}) } \right] \\
        & = \frac{1}{2} \left[ (\overline{\varphi} u_{\varepsilon}) \left( \begin{matrix}
             -\mathrm{i} (\partial_x \varphi) u_{\varepsilon}  - \pi \beta  f |\varphi|^2 \varphi u_{\varepsilon}  \\
              -\mathrm{i} \varphi \partial_y  u_{\varepsilon}
        \end{matrix}\right)  +  c.c.   \right] \\
        & = \frac{1}{2} \left( \begin{matrix}
            - \mathrm{i} \overline{\varphi} \partial_x \varphi +\mathrm{i} \varphi \partial_x \overline{\varphi} - 2 \pi \beta  f |\varphi|^4 \\
            0
        \end{matrix} \right) u_{\varepsilon}^2.
    \end{split}
\end{equation}
Then we have
\begin{align*}
        &- 2\beta  \varphi \int_{\mathbb{R}} \bT_0 \star \bj  \, u_{\varepsilon}^2 \mathrm{d} y = \pi \beta \varphi \\
        &\int_{\mathbb{R}}  \int_{\mathbb{R}} \operatorname{sgn}(y-y') \left( - \mathrm{i} \overline{\varphi} \partial_x \varphi + c.c. - 2 \pi \beta  f(y') |\varphi|^4 \right) \\ &\quad\quad u_{\varepsilon}^2(y') \mathrm{d} y'  u_{\varepsilon}^2 (y)  \mathrm{d} y \\
        & = \pi \beta \varphi \int_{\mathbb{R}}\left( \mathrm{i} \overline{\varphi} \partial_x \varphi  +c.c. + 2 \pi \beta  f(y') |\varphi|^4 \right) f(y') u_{\varepsilon}^2 (y')  \mathrm{d} y' \\
        & = \frac{2}{3}\pi^2 \beta^2 |\varphi|^4 \varphi
   \end{align*}
using~\eqref{eq:f}-\eqref{eq:f3} again. This is indeed~\eqref{eq:last claim}, and there remains to insert in~\eqref{eq:preconc} to deduce~\eqref{eq:1DNLS}. We next discuss the effective theory~\eqref{eq:1DNLS} per se.

\bigskip 

\noindent\textbf{Approximate ground states of 1D NLS.} Our convention is that states are $L^2$-normalized, $\int_{\R} |\varphi |^2 = 1$. In the aforementioned cold atoms context this means that we have $\beta = \alpha N$ with $N$ the particle number, $\alpha$ the (artificial) magnetic flux per particle, and 
$$\gt = g \varepsilon ^{-1/2} \int_\R |u_1|^4 \propto N a \varepsilon ^{-1/2}$$
with $a$ the scattering length of the pair-interaction potential. We first discuss the Thomas-Fermi (TF) approximation to in-trap ground states of~\eqref{eq:1Dener}. Namely we solve the variational equation for energy minimizers by neglecting the kinetic energy/derivative term, leading to 
\begin{equation}\label{eq:TF}
\pi^2 \beta^2 \varrho ^{2} + \gt \varrho + x^2 = \mu 
\end{equation}
in terms of the density $\varrho = |\varphi|^2$ and with a chemical potential $\mu$ ensuring normalization. For a solution living on a length scale $L$, a rough optimization of the energy as a function of $L$ shows that the effect of the quintic NLS term is dominant for $\beta^{3/2} \gg |\gt|$. We focus on this regime, which better shows the influence of the original flux attachment in the approximate density
\begin{equation}\label{eq:TFdens}
\varrho^{\rm TF} = \frac{1}{\pi |\beta|} \left((2|\beta|)^{1/2} -x^2\right)^{1/2} 
\end{equation}
where we have neglected the correction due to pair-interactions. In a quadratic trap the density's width is thus of order $|\beta|^{1/4}$ and the density has a distinctive semi-circle-law shape. We refer to~\cite{FeiChe-24} for more discussion.

\bigskip 

\noindent\textbf{Soliton solutions of 1D NLS in the focusing case.} An interesting possibility arises if the pair interactions are made attractive, $\gt <0$, and the gas is released from the trap. Indeed, the free 1D NLS
\begin{equation}\label{eq:freeNLS}
\im \partial_t \varphi = - \partial_x ^2 \varphi + \pi ^2 \beta ^2 |\varphi|^4 \varphi + \gt |\varphi|^2 \varphi  
\end{equation}
has in this case exact traveling waves solutions
$$ \varphi (x,t) \propto e^{\im ( p t + q x)} u(x-ct).$$
More precisely, one finds~\cite{PusPusTom-79,CowEnsRanSan-86,GraGra-92,GenMalWei-16} solutions of the form 
\begin{multline}\label{eq:soliton}
\varphi_{\rm sol} (x,t) = \\\frac{C}{A^{1/2}} e^{\im \left(\frac{2\omega t}{A^2} + \frac{kx}{A}\right)}\left(1+D\cosh\left(2C \left(\frac{x}{A}-\frac{2tk}{A^2}\right)\right) \right)^{-1/2} 
\end{multline}
with 
\begin{equation}\label{eq:solpara}
A= \frac{4}{|\gt|},\quad  C^2 = 2 \omega + k^2, \quad D^2 = 1 - \frac{\pi}{3} \beta^2 C^2.
\end{equation}
The solution exists provided that $C,D$ can be chosen non-negative. Imposing the mass-constraint $\int_{\R} |\varphi_{\rm sol}|^2 = 1$ by using that~\cite[Equation~(5)]{BirMal-08} 
$$ \int_{R} \frac{1}{1+D\cosh\left( 2C y\right)} \mathrm{d} y = \frac{C}{2\sqrt{1-D^2}}\log \frac{1+\sqrt{1-D^2}}{1-\sqrt{1-D^2}}$$
we find 
$$ C = \sqrt{\frac{3}{\pi}} |\beta|^{-1} \frac{\exp\left(2\sqrt{\frac{\pi}{3}}|\beta|\right)-1}{\exp\left(2\sqrt{\frac{\pi}{3}}|\beta|\right)+1},$$
which indeed ensures that $C,D >0$. Hence the soliton's width is of the order 
$$ L_{\rm sol} = \frac{A}{C} \sim 4\sqrt{\frac{3}{\pi}}\frac{|\beta|}{|\gt|}$$
for large $|\beta|$, and 
$$ L_{\rm sol} \sim \frac{4}{|\gt|}$$ 
for small $|\beta|$, where one recovers the usual soliton of the focusing cubic NLS. The former regime is of course more interesting in showing the 1D remnant of the original flux attachment. In particular we have 
$$|D| \underset{|\beta|\to\infty}{\ll} 1$$
so that the density $\left|\varphi_{\rm sol}\right|^2$ of the soliton~\eqref{eq:soliton} has a distinctively flat bulk profile before quickly dropping to $0$ at the edge. This contrasts with the soliton profile of the purely cubic focusing NLS, and is reminiscent of a quantum droplet phase found in ultradilute 1D Bose gases~\cite{PetAst-16,AstMal-18}. In the latter case however the effective equation incorporating beyond mean-field effects is a cubic-defocusing/quadratic-focusing NLS, whose soliton profile differs from the above. 

\bigskip 

\noindent\textbf{Concluding remarks.} We have derived a 1D NLS equation~\eqref{eq:1DNLS} as the dimensional reduction of the mean-field description of a macroscopic system of 2D anyons~\eqref{eq:2DPDE}, i.e. the Chern-Simons-Schr\"odinger system. In previous contributions~\cite{RouYan-23a,RouYan-23b} we have shown that the many-body problem for anyons converges to a different theory, namely the Tonks-Girardeau Bose gas, in the 1D limit at \emph{fixed particle number}. Thus the 1D and macroscopic/mean-field/almost bosonic limits do not commute, and this begs the question of what limit should be taken first in a given practical situation. In any case, none of the proposed 1D anyon models we are aware of is the dimensional reduction of the 2D anyon model in the sense we investigated here. Noticeably, a quintic NLS theory similar to ours has been proposed as an effective model for a 1D Calogero-Sutherland gas~\cite{FeiChe-24}, which is one of the proposed models for 1D fractional statistics.

The situation of a tight waveguide acting on the Schr\"odinger field that we considered takes inspiration from the typical cold atoms experimental set-up. External potentials such as magneto-optics traps acting on alkali gases can indeed be manipulated efficiently. This set-up for dimensional reduction differs markedly from other approaches~\cite{HanLeiMyr-92,AglGriJacPiSem-96,ChiEtalCel-22,ValOhb-24} that lead to other models, some of which have been experimentally implemented~\cite{FroEtalTar-22,DhaEtalNag-24,KwaEtalGre-24}. 

The main difference in our approach is that we directly reduce to 1D the full Chern-Simons-Schr\"odinger system, instead of first reducing the Chern-Simons action, and then coupling the obtained model to a 1D matter field as in~\cite{ValOhb-24}. We feel this is more natural for the 1D confinement of a 2D Bose-Einstein condensate coupled to artificial magnetic flux, but probably less relevant for boundary theories of FQH fluids.  

In essence, the difference between our effective model and those of ~\cite{AglGriJacPiSem-96,ChiEtalCel-22,ValOhb-24,Rabello-96,Kundu-99,BatGuaOel-06,GreSan-15,TanEggPel-15,Bon-etal-Pos-21,KwaEtalGre-24} stems from the fact that we have a magnetic \emph{field} proportional to matter density~\cite{ValWesOhb-20}  as opposed to a magnetic \emph{vector potential} proportional to matter density~\cite{EdmEtalOhb-13}. During the 1D reduction, it then becomes favorable to develop a specific phase factor~\eqref{eq:changegauge}, which leads to an effective purely local theory despite the initial long-range 2D magnetic interactions.

Our discussion can provide some experimental tests of artificial magnetic flux attachment. Loading a gas coupled to a density-dependent gauge field in a tight waveguide would allow to compare data to the well-known properties (some of which we recalled above) of the 1D NLS~\eqref{eq:1DNLS}, where the magnetic flux coupling constant $\beta$ appears explicitly. It particular it sets the width a semi-circle density profile in a harmonic trap. If the two-body coupling constant $g$ of the gas can be made negative (attractive interaction),~\eqref{eq:1DNLS} has soliton solutions upon removing the $x^2$ trapping potential in the loose direction. They differ from those observed~\cite{FroEtalTar-22} with a density-dependent vector potential. They are not chiral, and have a different, markedly flat, density profile.

\bigskip
\noindent\textbf{Acknowledgments.} We thank Blagoje Oblak, Douglas Lundholm, Patrik \"Ohberg, Gerard Valent\'{\i}-Rojas, Alessio Celi and Claudio Iacovelli for stimulating discussions.

\newpage

%

\end{document}